# CRITICALITY FOR BRANCHING PROCESSES IN RANDOM ENVIRONMENT[1]

By V. I. Afanasyev, J. Geiger, G. Kersting and V. A. Vatutin

*Steklov Institute, University Kaiserslautern,
University Frankfurt/Main and Steklov Institute*

We study branching processes in an i.i.d. random environment, where the associated random walk is of the oscillating type. This class of processes generalizes the classical notion of criticality. The main properties of such branching processes are developed under a general assumption, known as Spitzer's condition in fluctuation theory of random walks, and some additional moment condition. We determine the exact asymptotic behavior of the survival probability and prove conditional functional limit theorems for the generation size process and the associated random walk. The results rely on a stimulating interplay between branching process theory and fluctuation theory of random walks.

**1. Introduction and main results.** In this paper fundamental properties of branching processes in a critical random environment are developed. In such a process individuals reproduce independently of each other according to random offspring distributions which vary from one generation to the other. To give a formal definition let $\Delta$ be the space of probability measures on $\mathbb{N}_0 := \{0, 1, \dots\}$. Equipped with the metric of total variation, $\Delta$ becomes a Polish space. Let $Q$ be a random variable taking values in $\Delta$. Then, an infinite sequence $\Pi = (Q_1, Q_2, \dots)$ of i.i.d. copies of $Q$ is said to form a *random environment*. A sequence of $\mathbb{N}_0$-valued random variables $Z_0, Z_1, \dots$ is called a *branching process in the random environment* $\Pi$, if $Z_0$ is independent of $\Pi$ and given $\Pi$ the process $Z = (Z_0, Z_1, \dots)$ is a Markov chain with

$$(1.1) \qquad \mathcal{L}(Z_n | Z_{n-1} = z, \Pi = (q_1, q_2, \dots)) = \mathcal{L}(\xi_1 + \cdots + \xi_z)$$

Received November 2003; revised April 2004.
[1]Supported in part by the German Research Foundation (DFG) and the Russian Foundation of Basic Research Grant 436 Rus 113/683 and in part by Grant 02-01-00266.

for every $n \geq 1, z \in \mathbb{N}_0$ and $q_1, q_2, \ldots \in \Delta$, where $\xi_1, \xi_2, \ldots$ are i.i.d. random variables with distribution $q_n$. In the language of branching processes $Z_n$ is the $n$th generation size of the population and $Q_n$ is the distribution of the number of children of an individual at generation $n-1$.

We will denote the corresponding probability measure on the underlying probability space by $\mathbf{P}$. (If we refer to other probability spaces, then we use notation $\mathbb{P}$ and $\mathbb{E}$ for the respective probability measures and expectations.) Property (1.1) can be equivalently expressed as

$$(1.2)\ \mathbf{P}\{(Z_1, \ldots, Z_n) \in B | Z_0 = z_0, \Pi\} = \kappa_{n,z_0}(Q_1, \ldots, Q_n; B), \qquad \mathbf{P}\text{-a.s.},$$

where $B \subset \mathbb{N}_0^n$ and $\kappa_{n,z_0}$ is the kernel

$$\kappa_{n,z_0}(q_1, \ldots, q_n; B) := \sum_{(z_1, \ldots, z_n) \in B} q_1^{*z_0}(\{z_1\}) \cdots q_n^{*z_{n-1}}(\{z_n\}),$$

built from the $z$-fold convolutions $q_i^{*z}$ of the $q_i$. In the theorems below we assume $Z_0 = 1$ a.s. for convenience. Sometimes it will be necessary to allow other values $z$ for $Z_0$. Then, as usual we write $\mathbf{P}_z\{\cdot\}$ and $\mathbf{E}_z[\cdot]$ for the corresponding probabilities and expectations. For further details and background we refer the reader to Athreya and Karlin (1971), Athreya and Ney (1972) and Smith and Wilkinson (1969).

As it turns out the properties of $Z$ are first of all determined by its associated random walk $S = (S_0, S_1, \ldots)$. This random walk has initial state $S_0 = 0$ and increments $X_n = S_n - S_{n-1}$, $n \geq 1$, defined as

$$X_n := \log m(Q_n),$$

which are i.i.d. copies of the logarithmic mean offspring number $X := \log m(Q)$ with

$$m(Q) := \sum_{y=0}^{\infty} y Q(\{y\}).$$

We assume that $X$ is a.s. finite. In view of (1.1) the conditional expectation of $Z_n$, given the environment $\Pi$,

$$\mu_n := \mathbf{E}[Z_n | Z_0, \Pi]$$

can be expressed by means of $S$ as

$$\mu_n = Z_0 e^{S_n}, \qquad \mathbf{P}\text{-a.s.}$$

According to fluctuation theory of random walks [cf. Chapter XII in Feller (1971)], one may distinguish three different types of branching processes in a random environment. First, $S$ can be a random walk with positive drift, which means that $\lim_{n \to \infty} S_n = \infty$ a.s. In this case one has $\mu_n \to \infty$ a.s., provided $Z_0 \geq 1$, and $Z$ is called a *supercritical* branching process. Second, $S$ can have negative drift, that is, $\lim_{n\to\infty} S_n = -\infty$ a.s. Then $\mu_n \to 0$ a.s. and $Z$ is called *subcritical*. Finally, $S$ may be an oscillating random walk, meaning that $\limsup_{n\to\infty} S_n = \infty$ a.s. and, at the same time, $\liminf_{n\to\infty} S_n = -\infty$ a.s., which implies $\limsup_{n\to\infty} \mu_n = \infty$ a.s. and $\liminf_{n\to\infty} \mu_n = 0$ a.s. Then we call $Z$ a *critical* branching process. Our classification extends the classical distinction of branching processes in random environment introduced in Athreya and Karlin (1971) and Smith and Wilkinson (1969). There it is assumed that the random walk has finite mean. In this case the branching process $Z$ is supercritical, subcritical or critical, depending on whether $\mathbf{E}X > 0$, $\mathbf{E}X < 0$ or $\mathbf{E}X = 0$. Here we do not require that the expectation of $X$ exists.

In this paper we focus on the critical case, where the population eventually becomes extinct with probability 1. Indeed, note that the estimate

$$(1.3) \quad \begin{aligned} \mathbf{P}\{Z_n > 0 | Z_0, \Pi\} &= \min_{m \leq n} \mathbf{P}\{Z_m > 0 | Z_0, \Pi\} \\ &\leq \min_{m \leq n} \mu_m = Z_0 \exp\left(\min_{m \leq n} S_m\right) \end{aligned}$$

implies $\mathbf{P}\{Z_n > 0 | Z_0, \Pi\} \to 0$ a.s. in critical (as well as subcritical) cases and, consequently,

$$\mathbf{P}\{Z_n > 0\} \to 0 \qquad \text{as } n \to \infty.$$

A main task in the investigation of critical branching processes consists in determining the asymptotic probability of the event $\{Z_n > 0\}$ of nonextinction at generation $n$ and the asymptotic behavior of $Z$ on this event. To this end, we impose an assumption on the random walk $S$, which is known as *Spitzer's condition*. This condition says that the expected proportion of time, which the random walk spends within the positive real half line up to time $n$, stabilizes as $n \to \infty$ at some value other than 0 or 1.

ASSUMPTION A1. There exists a number $0 < \rho < 1$ such that

$$\frac{1}{n} \sum_{m=1}^{n} \mathbf{P}\{S_m > 0\} \to \rho \qquad \text{as } n \to \infty.$$

This condition plays an important role in fluctuation theory of random walks. The summands $\mathbf{P}\{S_m > 0\}$ may likewise be replaced by $\mathbf{P}\{S_m \geq 0\}$, since $\sum_{m=1}^{n} \mathbf{P}\{S_m = 0\} = o(n)$ for every nondegenerate random walk $S$ [cf. XII.9(c) in Feller (1971)]. In fact, Doney (1995) proved that Assumption A1 is equivalent to the convergence of $\mathbf{P}\{S_n > 0\}$ to $\rho$. It is well known that any random walk satisfying Assumption A1 is of the oscillating type [see, e.g., Section XII.7 in Feller (1971)]. We note that Assumptions A1 covers nondegenerate random walks with zero mean and finite variance increments, as well as all nondegenerate symmetric random walks. In these cases $\rho = 1/2$. Other examples are provided by random walks in the domain of attraction of some stable law, see Assumption B1.

Our second assumption on the environment concerns the standardized truncated second moment of $Q$,

$$\zeta(a) := \sum_{y=a}^{\infty} y^2 Q(\{y\})/m(Q)^2, \qquad a \in \mathbb{N}_0.$$

To formulate the assumption let us introduce the renewal function

(1.4) $$v(x) := \begin{cases} 1 + \sum_{i=1}^{\infty} \mathbf{P}\{S_{\gamma_i} \geq -x\}, & \text{if } x \geq 0, \\ 0, & \text{else,} \end{cases}$$

where $0 =: \gamma_0 < \gamma_1 < \cdots$ are the strict descending ladder epochs of $S$,

$$\gamma_i := \min\{n > \gamma_{i-1} : S_n < S_{\gamma_{i-1}}\}, \qquad i \geq 1.$$

The fundamental property of $v$ is

(1.5) $$\mathbf{E}v(x + X) = v(x), \qquad x \geq 0,$$

which holds for any oscillating random walk [cf. Bertoin and Doney (1994) and Kozlov (1976)].

ASSUMPTION A2. For some $\epsilon > 0$ and some $a \in \mathbb{N}_0$,

$$\mathbf{E}(\log^+ \zeta(a))^{1/\rho + \epsilon} < \infty \quad \text{and} \quad \mathbf{E}[v(X)(\log^+ \zeta(a))^{1+\epsilon}] < \infty.$$

EXAMPLES. Here are some instances where this assumption is fulfilled:

1. If the random offspring distribution $Q$ has uniformly bounded support, that is, if $\mathbf{P}\{Q(\{0, 1, \ldots, a^*\}) = 1\} = 1$ for some $a^*$, then $\zeta(a) = 0$ $\mathbf{P}$-a.s. for all $a > a^*$. In this case Assumption A2 is redundant and we merely require the random walk $S$ to satisfy Spitzer's condition. In particular, the results to follow hold for any binary branching process in a random environment (where individuals have either two children or none), which satisfies Assumption A1.
2. In view of relation (1.5), we have $\mathbf{E}v(X) = v(0) < \infty$. Therefore, Assumption A2 is satisfied, if $\zeta(a)$ is a.s. bounded from above for some $a$. We note that this is the case if the value of $Q$ is a.s. a Poisson distribution or a.s. a geometric distribution on $\mathbb{N}_0$ (with random expectations). This follows from the estimate

$$\frac{\zeta(2)}{2} \leq \eta := \sum_{y=0}^{\infty} y(y-1)Q(\{y\})/m(Q)^2$$

and the observation that for Poisson distributions $\eta = 1$ a.s. and for geometric distributions on $\mathbb{N}_0$, one has $\eta = 2$ a.s.

3. The renewal function $v(x)$ always satisfies $v(x) = O(x)$ as $x \to \infty$ and $v(x) = 0$ for $x < 0$. Therefore, as follows from Hölder's inequality, Assumption A2 is entailed by

$$\mathbf{E}(X^+)^p < \infty \quad \text{and} \quad \mathbf{E}(\log^+ \zeta(a))^q < \infty$$

for some $p > 1$ and $q > \max(\rho^{-1}, p(p-1)^{-1})$.

If $X$ has regular tails we can replace Assumptions A1 and A2 by the following alternative set of conditions.

ASSUMPTION B1. The distribution of $X$ belongs without centering to the domain of attraction of some stable law $\lambda$ with index $\alpha \in (0, 2]$. The limit law $\lambda$ is not a one-sided stable law, that is, $0 < \lambda(\mathbb{R}^+) < 1$.

Here it is assumed that there are numbers $c_n > 0$ such that $c_n S_n$ converges in distribution to $\lambda$ and, consequently, $\mathbf{P}(S_n > 0) \to \rho := \lambda(\mathbb{R}^+)$. In particular, Assumption B1 implies Assumption A1.

The gain of the stronger regularity condition B1 is that we can further relax the integrability condition A2.

ASSUMPTION B2. For some $\epsilon > 0$ and some $a \in \mathbb{N}_0$,

$$\mathbf{E}(\log^+ \zeta(a))^{\alpha + \epsilon} < \infty.$$

We note that Assumption A2 is indeed stronger than Assumption B2 since

(1.6) $$\rho \leq \alpha^{-1}.$$

For further explanations the reader may consult Dyakonova, Geiger and Vatutin (2004) or Chapter 8.9.2 in the monograph of Bingham, Goldie and Teugels (1987).

We now come to the main results of the paper. All our limit theorems are under the law $\mathbf{P}$, which is what is called the annealed approach. The first theorem describes the asymptotic behavior of the nonextinction probability at generation $n$.

THEOREM 1.1. *Assume Assumptions* A1 *and* A2 *or* B1 *and* B2. *Then there exists a positive finite number $\theta$ such that*

(1.7) $$\mathbf{P}\{Z_n > 0\} \sim \theta \mathbf{P}\{\min(S_1, \ldots, S_n) \geq 0\} \qquad \text{as } n \to \infty.$$

66     V. I. AFANASYEV, J. GEIGER, G. KERSTING AND V. A. VATUTIN

This theorem gives first evidence for our claim that the asymptotic behavior of $Z$ is primarily determined by the random walk $S$, since only the constant $\theta$ depends on the fine structure of the random environment. The asymptotics (1.7) reflects the following fact: If $\min_{m \leq n} S_m$ is low, then the probability of nonextinction at $n$ is very small as follows from (1.3). In fact, it turns out that on the event $\{Z_n > 0\}$, the value of $\min_{m \leq n} S_m$ is only of constant order. A detailed description of this phenomenon is given in Theorem 1.4.

Since the asymptotic behavior of the probability on the right-hand side of (1.7) is well known under Assumption A1 (see Lemma 2.1), we obtain the following corollary.

COROLLARY 1.2. *Assume Assumptions* A1 *and* A2 *or* B1 *and* B2. *Then,*

$$\mathbf{P}\{Z_n > 0\} \sim \theta n^{-(1-\rho)} l(n) \qquad as \ n \to \infty,$$

*where $l(1), l(2), \ldots$ is a sequence varying slowly at infinity.*

A probabilistic representation of $\theta$ is contained in (4.10) and explained in the remark thereafter. For a representation of the function $l$, see Lemma 2.1.

The next theorem shows that conditioned on the event $\{Z_n > 0\}$, the generation size process $Z_0, Z_1, \ldots, Z_n$ exhibits "supercritical" behavior. Supercritical branching processes (whether classical or in a random environment) obey the growth law $Z_n / \mu_n \to W$ a.s., where $W$ is some typically nondegenerate random variable. However, in our situation this kind of behavior can no longer be formulated as a statement on a.s. convergence, since the conditional distribution of the environment $\Pi$, given $\{Z_n > 0\}$, changes with $n$.

Instead, let us, for integers $0 \leq r \leq n$, consider the rescaled generation size process $X^{r,n} = (X_t^{r,n})_{0 \leq t \leq 1}$, given by

$$(1.8) \qquad X_t^{r,n} := \frac{Z_{r + \lfloor (n-r)t \rfloor}}{\mu_{r + \lfloor (n-r)t \rfloor}}, \qquad 0 \leq t \leq 1.$$

THEOREM 1.3. *Assume Assumptions* A1 *and* A2 *or* B1 *and* B2. *Let $r_1, r_2, \ldots$ be a sequence of positive integers such that $r_n \leq n$ and $r_n \to \infty$. Then,*

$$\mathcal{L}(X^{r_n,n} | Z_n > 0) \implies \mathcal{L}((W_t)_{0 \leq t \leq 1}) \qquad as \ n \to \infty,$$

*where the limiting process is a stochastic process with a.s. constant paths, that is, $\mathbb{P}\{W_t = W \text{ for all } t \in [0,1]\} = 1$ for some random variable $W$. Furthermore,*

$$\mathbb{P}\{0 < W < \infty\} = 1.$$

Here, $\Longrightarrow$ denotes weak convergence w.r.t. the Skorokhod topology in the space $D[0,1]$ of càdlàg functions on the unit interval. Again, the growth of $Z$ is in the first place determined by the random walk [namely, the sequence $(\mu_n)_{n \geq 0}$]. The fine structure of the random environment is reflected only in the distribution of $W$.

Thus, first of all properties of the random walk $S$ are important for the behavior of $Z$. However, one also has to take into account that the random walk changes its properties drastically, when conditioned on the event $\{Z_n > 0\}$. The next theorem illustrates this fact. Let $\tau_n$ be the first moment, when the minimum of $S_0, \ldots, S_n$ is attained

$$(1.9) \qquad \tau_n := \min\{i \leq n | S_i = \min(S_0, \ldots, S_n)\}, \qquad n \geq 0.$$

THEOREM 1.4. *Assume Assumptions* A1 *and* A2 *or* B1 *and* B2. *Then, as* $n \to \infty$,

$$\mathcal{L}((\tau_n, \min(S_0, \ldots, S_n)) | Z_n > 0)$$

*converges weakly to some probability measure on* $\mathbb{N}_0 \times \mathbb{R}_0^-$.

For a more detailed description of the conditioned random walk we confine ourselves to the situation given in Assumption B1.

THEOREM 1.5. *Assume Assumptions* B1 *and* B2. *Then there exists a slowly varying sequence* $\ell(1), \ell(2), \ldots$ *such that*

$$\mathcal{L}((n^{-1/\alpha} \ell(n) S_{\lfloor nt \rfloor})_{0 \leq t \leq 1} | Z_n > 0) \quad \Longrightarrow \quad \mathcal{L}(L^+) \qquad as \ n \to \infty,$$

*where* $L^+$ *denotes the meander of a strictly stable process with index* $\alpha$.

Shortly speaking, the meander $L^+ = (L_t^+)_{0 \leq t \leq 1}$ is a strictly stable Lévy process conditioned to stay positive on the time interval $(0, 1]$ [for details see Doney (1985) and Durrett (1978)].

In view of Theorem 1.3, the last theorem is equivalent to the following result.

COROLLARY 1.6. *Assume Assumptions* B1 *and* B2. *Then,*

$$\mathcal{L}((n^{-1/\alpha} \ell(n) \log Z_{\lfloor nt \rfloor})_{0 \leq t \leq 1} | Z_n > 0) \quad \Longrightarrow \quad \mathcal{L}(L^+) \qquad as \ n \to \infty$$

*for some slowly varying sequence* $\ell(1), \ell(2), \ldots$.

Starting from the seminal paper of Kozlov (1976), the topic of branching processes in a critical random environment has gone through quite a development. For a fairly long time research was restricted to the special case of random offspring distributions with a linear fractional generating function (including geometric distributions) and to random walks with zero mean, finite variance increments. Under these restrictions, fairly explicit (albeit tedious) calculations of certain Laplace transforms are feasible, which then allow the proof of (most of) the results above [cf. Afanasyev (1993, 1997)]. In recent years the assumption of linear fractional offspring distributions could be dropped [see Afanasyev (2001), Geiger and Kersting (2000), Kozlov (1995) and Vatutin (2002)], and first steps to overcome the assumption of a finite variance random walk were taken [see Dyakonova, Geiger and Vatutin (2004) and Vatutin and Dyakonova (2003)].

Yet the significance of Spitzer's condition as a suitable regularity condition for branching processes in an i.i.d. random environment has not been recognized so far. The use of Laplace transforms and generating functions is still indispensable for our purposes (see Section 3), however, in our general situation it is to be supported by other devices. In particular, we point out to the change of measure, which is discussed in Section 2. It enables us to make use of Tanaka's path decomposition for conditioned random walks [Tanaka (1989)]. This decomposition turns out to be an essential tool in the proof of Theorem 1.3. In particular, it is used to establish the fact that the mean gives the right growth of the population up to a random factor, which was an open problem [see Afanasyev (2001)]. Tanaka's decomposition also allows to substantially weaken the required moment conditions to Assumption A2, respectively, Assumption B2.

**2. Auxiliary results for random walks.** The proofs of our theorems rely strongly on various results from the theory of random walks. In this section we collect these results.

2.1. *Results from fluctuation theory of random walks.* The minima

$$L_n := \min(S_1, \ldots, S_n), \qquad n \geq 1,$$

play an important role in the fluctuation theory of random walks. Recall the definition of the function $v(x)$ introduced in (1.4).

LEMMA 2.1. *Assume Assumption* A1. *Then, for every* $x \geq 0$,

(2.1) $\qquad \mathbf{P}\{L_n \geq -x\} \sim v(x) n^{-(1-\rho)} l(n) \qquad \text{as } n \to \infty,$

*where the slowly varying function $l$ is given by $l(n) := h(1 - n^{-1})/\Gamma(\rho)$ with*

$$h(s) := \exp\biggl(\sum_{n=1}^{\infty} \frac{s^n}{n}(\mathbf{P}\{S_n \geq 0\} - \rho)\biggr), \qquad 0 \leq s < 1.$$

*Furthermore, there exists a constant $0 < c_1 < \infty$ such that for all $x \geq 0$ and $n \in \mathbb{N}$,*

(2.2) $\qquad \mathbf{P}\{L_n \geq -x\} \leq c_1 v(x) n^{-(1-\rho)} l(n).$

*Similarly, for all $n$ and some $0 < c_2 < \infty$, we have*

(2.3) $\qquad \mathbf{P}\{\max(S_1, \ldots, S_n) \leq 0\} \leq c_2 n^{-\rho} l(n)^{-1}.$

PROOF. For the asymptotics (2.1) apply Theorem 8.9.12 in Bingham, Goldie and Teugels (1987) to the random walk $\overline{S} := -S$ and note that there $\rho$ has to be replaced by $\bar{\rho} := 1 - \rho$.

For the second claim [which has been established already by Kozlov (1976) for finite variance random walks] we use the inequality

$$\sum_{m=0}^{\infty} s^m \mathbf{P}\{L_m \geq -x\} \leq v(x) \exp\left(\sum_{n=1}^{\infty} \frac{s^n}{n} \mathbf{P}\{S_n \geq 0\}\right)$$
$$= v(x)(1-s)^{-\rho} h(s)$$

following from Lemma 8.9.11 and a formula contained in the proof of Theorem 8.9.12 in Bingham, Goldie and Teugels (1987). Since $\mathbf{P}\{L_n \geq -x\}$ is nonincreasing with $n$, it follows that

$$\frac{n}{2}\left(1 - \frac{1}{n}\right)^n \mathbf{P}\{L_n \geq -x\} \leq \sum_{n/2 \leq m \leq n} \left(1 - \frac{1}{n}\right)^m \mathbf{P}\{L_m \geq -x\}$$
$$\leq v(x) n^\rho h(1 - n^{-1}),$$

which implies the bound (2.2). Finally, (2.3) follows by applying (2.2) to the random walk $\overline{S}$. Then $\rho$ has to be replaced by $\bar{\rho} = 1 - \rho$, and $h(s)$ by

$$\bar{h}(s) := \exp\left(\sum_{n=1}^{\infty} \frac{s^n}{n}(\mathbf{P}\{S_n \leq 0\} - \bar{\rho})\right) = \exp\left(\sum_{n=1}^{\infty} \frac{s^n}{n}(\rho - \mathbf{P}\{S_n > 0\})\right).$$

As $s \uparrow 1$, we obtain

$$h(s)\bar{h}(s) = \exp\left(\sum_{n=1}^{\infty} \frac{s^n}{n} \mathbf{P}\{S_n = 0\}\right) \to \exp\left(\sum_{n=1}^{\infty} \frac{1}{n} \mathbf{P}\{S_n = 0\}\right) =: \gamma.$$

Since $\gamma$ is positive and finite [see XII.9(c) in Feller (1971)], we have $l(n)\bar{l}(n) \sim \gamma/(\Gamma(\rho)\Gamma(1-\rho))$, and the upper bound (2.3) follows. $\square$

Next we study the random time $\tau_n$ defined in (1.9). The following technical lemma will be used at various places.

LEMMA 2.2. *Let $u(x)$, $x \geq 0$, be a nonnegative, nonincreasing function with $\int_0^\infty u(x)\,dx < \infty$. Then, under Assumption* A1, *for every $\epsilon > 0$, there exists a positive integer $l$ such that for all $n \geq l$,*

$$\sum_{k=l}^{n} \mathbf{E}[u(-S_k); \tau_k = k] \mathbf{P}\{L_{n-k} \geq 0\} \leq \epsilon \mathbf{P}\{L_n \geq 0\}.$$

PROOF. We first show

$$\sum_{k=0}^{\infty} \mathbf{E}[u(-S_k); \tau_k = k] < \infty. \qquad (2.4)$$

Let $x \geq 0$. Since [recall (1.4) and note that $\mathbf{P}\{S_{\gamma_0} \geq -x\} = 1$ for $x \geq 0$]

$$v(x) = \sum_{k=0}^{\infty} \sum_{i=0}^{\infty} \mathbf{P}\{S_k \geq -x, \gamma_i = k\}$$
$$= \sum_{k=0}^{\infty} \mathbf{P}\{-S_k \leq x, \tau_k = k\}, \qquad (2.5)$$

we have

$$\sum_{k=0}^{\infty} \mathbf{E}[u(-S_k); \tau_k = k] = \int_0^{\infty} u(x)\, dv(x).$$

In view of the fact that $v(x) = O(x)$ as $x \to \infty$, assertion (2.4) now follows from the integrability and monotonicity assumptions on $u$. Next we prove

$$\mathbf{E}[u(-S_n); \tau_n = n] = O(n^{-1}). \qquad (2.6)$$

Clearly, $\tau_n = n$ implies $S_n \leq 0$. Hence,

$$\mathbf{E}[u(-S_n); \tau_n = n]$$
$$\leq \sum_{k=0}^{\infty} u(2^k - 1)\mathbf{P}\{-(2^{k+1} - 1) < S_n \leq -(2^k - 1), \tau_n = n\}$$
$$\leq \sum_{k=0}^{\infty} u(2^k - 1)\mathbf{P}\{S_0, \ldots, S_{\lfloor n/2 \rfloor} \geq -2^{k+1},$$
$$S_m - S_{\lfloor (n+1)/2 \rfloor} \geq S_n - S_{\lfloor (n+1)/2 \rfloor}, \lfloor (n+1)/2 \rfloor \leq m < n\}$$
$$= \sum_{k=0}^{\infty} u(2^k - 1)\mathbf{P}\{L_{\lfloor n/2 \rfloor} \geq -2^{k+1}\}\mathbf{P}\{S_0, \ldots, S_{\lfloor n/2 \rfloor - 1} \geq S_{\lfloor n/2 \rfloor}\}.$$

By duality, $\mathbf{P}\{S_0, \ldots, S_{m-1} \geq S_m\} = \mathbf{P}\{\max(S_1, \ldots, S_m) \leq 0\}$. Therefore, using the upper bounds (2.2) and (2.3), we deduce

$$\mathbf{E}[u(-S_n); \tau_n = n]$$
$$\leq c_1 c_2 \lfloor n/2 \rfloor^{\rho-1} l(\lfloor n/2 \rfloor) \lfloor n/2 \rfloor^{-\rho} l(\lfloor n/2 \rfloor)^{-1} \sum_{k=0}^{\infty} u(2^k - 1) v(2^{k+1}).$$

Since $v(x) = O(x)$, the series is convergent by the assumption on $u$, and assertion (2.6) follows. Now observe that, by (2.6) and monotonicity of $L_j$, we have for any $0 < \delta < 1$,

$$\sum_{k=l}^{n} \mathbf{E}[u(-S_k); \tau_k = k] \mathbf{P}\{L_{n-k} \geq 0\}$$

$$\leq \mathbf{P}\{L_{\lfloor \delta n \rfloor} \geq 0\} \sum_{l \leq k \leq (1-\delta)n} \mathbf{E}[u(-S_k); \tau_k = k] + \frac{c}{(1-\delta)n} \sum_{j \leq \delta n} \mathbf{P}\{L_j \geq 0\},$$

where $c$ is some positive finite constant. By Lemma 2.1, $\mathbf{P}\{L_n \geq 0\}$ is regularly varying with exponent $-(1-\rho) \in (-1, 0)$. An application of Karamata's theorem [see, e.g., Theorem 1.5.11 in Bingham, Goldie and Teugels (1987)] gives

$$\sum_{j \leq \delta n} \mathbf{P}\{L_j \geq 0\} \sim \frac{\delta n}{\rho} \mathbf{P}\{L_{\lfloor \delta n \rfloor} \geq 0\} \sim \frac{\delta^\rho}{\rho} n \mathbf{P}\{L_n \geq 0\} \qquad \text{as } n \to \infty.$$

Consequently,

$$\sum_{k=l}^{n} \mathbf{E}[u(-S_k); \tau_k = k] \mathbf{P}\{L_{n-k} \geq 0\}$$

$$\leq c \mathbf{P}\{L_n \geq 0\} \left( \delta^{-(1-\rho)} \sum_{k=l}^{\infty} \mathbf{E}[u(-S_k); \tau_k = k] + \frac{\delta^\rho}{1-\delta} \right)$$

for sufficiently large $c$. By (2.4), the sum on the right-hand side above is finite. Hence, the claim of the lemma follows by a suitable choice of $\delta$ and $l$. □

As an application of Lemma 2.2 we generalize a functional limit theorem for random walks satisfying Assumption B1, which is due to Doney (1985) and Durrett (1978).

LEMMA 2.3. *Assume Assumption* B1 *and let* $x \geq 0$. *Then, there exists a slowly varying sequence* $\ell(1), \ell(2), \ldots$ *such that*

$$\mathcal{L}((n^{-1/\alpha} \ell(n) S_{\lfloor nt \rfloor})_{0 \leq t \leq 1} | L_n \geq -x) \quad \Longrightarrow \quad \mathcal{L}(L^+) \qquad \text{as } n \to \infty,$$

*where* $L^+$ *is the meander of a strictly stable Lévy process.*

PROOF. Doney and Durrett proved this theorem for $x = 0$. To treat the general case let us consider the processes $S^{k,n}$ and $\widetilde{S}^{k,n}, 0 \leq k \leq n$, given by

(2.7)
$$S_t^{k,n} := n^{-1/\alpha} \ell(n) S_{\lfloor nt \rfloor \wedge k},$$
$$\widetilde{S}_t^{k,n} := n^{-1/\alpha} \ell(n) (S_{\lfloor nt \rfloor} - S_{\lfloor nt \rfloor \wedge k}), \qquad 0 \leq t \leq 1.$$

Then $S^n := S^{k,n} + \widetilde{S}^{k,n}$ is the process under consideration. For $0 \leq k \leq n$, we define

$$(2.8) \qquad L_{k,n} := \min_{0 \leq j \leq n-k}(S_{k+j} - S_k).$$

Let $\phi$ be a bounded continuous function on the space $D[0,1]$ of càdlàg functions equipped with the Skorokhod metric. Since

$$(2.9) \qquad \{\tau_n = k\} = \{\tau_k = k\} \cap \{L_{k,n} \geq 0\},$$

a decomposition according to $\tau_n$ gives

$$(2.10) \quad \mathbf{E}[\phi(S^n); L_n \geq -x] = \sum_{k=0}^{n} \mathbf{E}[\phi(S^n); \tau_k = k, S_k \geq -x, L_{k,n} \geq 0].$$

Observe that, for $0 \leq k \leq n$, we have

$$(2.11) \quad \begin{aligned} & \mathbf{E}[\phi(S^n); \tau_k = k, S_k \geq -x, L_{k,n} \geq 0] \\ & = \mathbf{E}[\mathbf{E}[\phi(S^{k,n} + \widetilde{S}^{k,n}); L_{k,n} \geq 0 | X_1, \ldots, X_k]; \tau_k = k, S_k \geq -x]. \end{aligned}$$

For $k \geq 0$ fixed the result of Doney and Durrett implies that, given $L_{k,n} \geq 0$, the process $\widetilde{S}^{k,n}$ converges in distribution to the specified meander. Also, given $X_1, \ldots, X_k$, the process $S^{k,n}$ vanishes asymptotically $\mathbf{P}$-a.s. Hence, by independence and the dominated convergence theorem, we have

$$(2.12) \quad \begin{aligned} & \mathbf{E}[\phi(S^n); \tau_k = k, S_k \geq -x, L_{k,n} \geq 0] \\ & = \mathbf{P}\{L_n \geq 0\} \mathbf{P}\{\tau_k = k, S_k \geq -x\}(\mathbb{E}\phi(L^+) + o(1)) \end{aligned}$$

for every $k \geq 0$. Now let $\varepsilon > 0$. Taking $u = \mathbf{1}_{[0,x]}$ in Lemma 2.2 gives

$$(2.13) \quad \begin{aligned} & \sum_{k=l}^{n} \mathbf{E}[\phi(S^n); \tau_k = k, S_k \geq -x, L_{k,n} \geq 0] \\ & \leq \sup |\phi| \sum_{k=l}^{n} \mathbf{E}[u(-S_k); \tau_k = k] \mathbf{P}\{L_{n-k} \geq 0\} \leq \epsilon \mathbf{P}\{L_n \geq 0\}, \end{aligned}$$

if only $l$ is large enough. Combining formulas (2.10), (2.12) and (2.13) with (2.5) gives

$$\mathbf{E}[\phi(S^n); L_n \geq -x] = \mathbf{P}\{L_n \geq 0\} v(x)(\mathbb{E}\phi(L^+) + o(1)).$$

In particular, choosing $\phi \equiv 1$, we obtain the asymptotics for $\mathbf{P}\{L_n \geq -x\}$ and the claim of the lemma follows. $\square$

2.2. *A change of measure.* Following Geiger and Kersting (2000), it is helpful to consider, besides **P**, another probability measure $\mathbf{P}^+$. In order to define this measure let $\mathcal{F}_n, n \geq 0$, be the $\sigma$-field of events generated by the random variables $Q_1, \ldots, Q_n$ and $Z_0, \ldots, Z_n$. These $\sigma$-fields form a filtration $\mathcal{F}$.

LEMMA 2.4. *The random variables* $v(S_n)I_{\{L_n \geq 0\}}$, $n = 0, 1, \ldots$ *form a martingale with respect to* $\mathcal{F}$ *under* **P**.

PROOF. Let $B$ and $D$ be Borel sets in $\mathbb{N}_0^n$ and $\Delta^n$, respectively. Recall identities (1.2) and (1.5) and the fact that $v(x) = 0$ for $x < 0$. Conditioning first on the environment $\Pi$ and then on $\mathcal{F}_n$ and using the independence of $Q_1, Q_2, \ldots$, we obtain

$$\mathbf{E}[v(S_{n+1}); L_{n+1} \geq 0, Z_0 = z, (Q_1, \ldots, Q_n) \in D, (Z_1, \ldots, Z_n) \in B]$$
$$= \mathbf{E}[v(X_{n+1} + S_n)\kappa_{n,z}(Q_1, \ldots, Q_n; B);$$
(2.14) $\qquad L_n \geq 0, Z_0 = z, (Q_1, \ldots, Q_n) \in D]$
$$= \mathbf{E}[v(S_n)\kappa_{n,z}(Q_1, \ldots, Q_n; B); L_n \geq 0, Z_0 = z, (Q_1, \ldots, Q_n) \in D]$$
$$= \mathbf{E}[v(S_n); L_n \geq 0, Z_0 = z, (Q_1, \ldots, Q_n) \in D, (Z_1, \ldots, Z_n) \in B].$$

By definition of conditional expectation, (2.14) implies

$$\mathbf{E}[v(S_{n+1})I_{\{L_{n+1} \geq 0\}}|\mathcal{F}_n] = v(S_n)I_{\{L_n \geq 0\}}, \qquad \mathbf{P}\text{-a.s.},$$

which is the desired martingale property. $\square$

Taking into account $v(0) = 1$, we may introduce probability measures $\mathbf{P}_n^+$ on the $\sigma$-fields $\mathcal{F}_n$ by means of the densities

$$d\mathbf{P}_n^+ := v(S_n)I_{\{L_n \geq 0\}}\, d\mathbf{P}.$$

Because of the martingale property, the measures are consistent, that is, $\mathbf{P}_{n+1}^+|\mathcal{F}_n = \mathbf{P}_n^+$. Therefore (choosing a suitable underlying probability space), there exists a probability measure $\mathbf{P}^+$ on the $\sigma$-field $\mathcal{F}_\infty := \bigvee_n \mathcal{F}_n$ such that

(2.15) $\qquad \mathbf{P}^+|\mathcal{F}_n = \mathbf{P}_n^+, \qquad n \geq 0.$

We note that (2.15) can be rewritten as

(2.16) $\qquad \mathbf{E}^+ Y_n = \mathbf{E}[Y_n v(S_n); L_n \geq 0]$

for every $\mathcal{F}_n$-measurable nonnegative random variable $Y_n$. This change of measure is the well-known Doob $h$-transform from the theory of Markov processes. In particular, under $\mathbf{P}^+$, the process $S$ becomes a Markov chain with state space $\mathbb{R}_0^+$ and transition kernel

$$P^+(x; dy) := \frac{1}{v(x)}\mathbf{P}\{x + X \in dy\}v(y), \qquad x \geq 0.$$

In our context $\mathbf{P}^+$ arises from conditioning:

LEMMA 2.5. *Assume Assumption* A1. *For* $k \in \mathbb{N}$, *let* $Y_k$ *be a bounded real-valued* $\mathcal{F}_k$-*measurable random variable. Then, as* $n \to \infty$,

$$\mathbf{E}[Y_k | L_n \geq 0] \to \mathbf{E}^+ Y_k.$$

*More generally, let* $Y_1, Y_2, \ldots$ *be a uniformly bounded sequence of real-valued random variables adapted to the filtration* $\mathcal{F}$, *which converges* $\mathbf{P}^+$-*a.s. to some random variable* $Y_\infty$. *Then, as* $n \to \infty$,

$$\mathbf{E}[Y_n | L_n \geq 0] \to \mathbf{E}^+ Y_\infty.$$

PROOF. For $x \geq 0$, write $m_n(x) := \mathbf{P}\{L_n \geq -x\}$. Then, for $k \leq n$, conditioning on $\mathcal{F}_k$ gives

$$\mathbf{E}[Y_k | L_n \geq 0] = \mathbf{E}\left[Y_k \frac{m_{n-k}(S_k)}{m_n(0)}; L_k \geq 0\right].$$

The first claim now follows from the asymptotics (2.1) and (2.2), the dominated convergence theorem and relation (2.16). For the second claim let $\gamma > 1$. Using again (2.1), (2.2) and (2.16), we obtain, for $k \leq n$,

$$|\mathbf{E}[Y_n - Y_k | L_{\lfloor \gamma n \rfloor} \geq 0]| \leq \mathbf{E}\left[|Y_n - Y_k| \frac{m_{\lfloor (\gamma-1)n \rfloor}(S_n)}{m_{\lfloor \gamma n \rfloor}(0)}; L_n \geq 0\right]$$

$$\leq c \left(\frac{\gamma - 1}{\gamma}\right)^{-(1-\rho)} \mathbf{E}[|Y_n - Y_k| v(S_n); L_n \geq 0]$$

$$= c \left(\frac{\gamma - 1}{\gamma}\right)^{-(1-\rho)} \mathbf{E}^+ |Y_n - Y_k|,$$

where $c$ is some positive constant. Letting first $n \to \infty$ and then $k \to \infty$, the right-hand side vanishes by the dominated convergence theorem. Thus, using the first part of the lemma, we conclude

$$\mathbf{E}[Y_n; L_{\lfloor \gamma n \rfloor} \geq 0] = (\mathbf{E}^+ Y_\infty + o(1)) \mathbf{P}\{L_{\lfloor \gamma n \rfloor} \geq 0\}.$$

Consequently, for some $c > 0$,

$$|\mathbf{E}[Y_n; L_n \geq 0] - \mathbf{E}^+ Y_\infty \mathbf{P}\{L_n \geq 0\}|$$
$$\leq |\mathbf{E}[Y_n; L_{\lfloor \gamma n \rfloor} \geq 0] - \mathbf{E}^+ Y_\infty \mathbf{P}\{L_{\lfloor \gamma n \rfloor} \geq 0\}| + c \mathbf{P}\{L_n \geq 0, L_{\lfloor \gamma n \rfloor} < 0\}$$
$$\leq (o(1) + c(1 - \gamma^{-(1-\rho)})) \mathbf{P}\{L_n \geq 0\},$$

where for the last inequality we also used (2.1) again. Since $\gamma$ may be chosen arbitrarily close to 1, we have

$$\mathbf{E}[Y_n; L_n \geq 0] - \mathbf{E}^+ Y_\infty \mathbf{P}\{L_n \geq 0\} = o(\mathbf{P}\{L_n \geq 0\}),$$

which is the second claim of the lemma. □

The change of measure has a natural interpretation, as is known from Bertoin and Doney (1994) and others: Under $\mathbf{P}^+$, the chain $S$ can be viewed as a random walk conditioned to never hit the strictly negative half line. Then $S$ gains an important renewal property, which is a consequence of the *Tanaka decomposition* for oscillating random walks. We recall only those aspects of the decomposition, which will be needed in the sequel and which, in our context, have to be extended to the entire environment. The original decomposition in Tanaka (1989) is not fully suitable for our purposes, since it concerns random works conditioned to never leave the strictly (!) positive real half line, meaning that it is based on a somewhat different harmonic function than $v(x)$. For these reasons, as well as for the readers convenience, we briefly recall the decomposition and its proof.

Let $\nu \geq 1$ be the time of the first *prospective minimal value* of $S$, that is, a minimal value with respect to the future development of the walk,

$$\nu := \min\{m \geq 1 : S_{m+i} \geq S_m \text{ for all } i \geq 0\}.$$

Moreover, let $\iota \geq 1$ be the first weak ascending ladder epoch of $S$,

$$\iota := \min\{m \geq 1 : S_m \geq 0\}.$$

We denote

$$\widetilde{Q}_n := Q_{\nu+n} \quad \text{and} \quad \widetilde{S}_n := S_{\nu+n} - S_\nu, \qquad n \geq 1.$$

LEMMA 2.6. *Suppose that $\iota < \infty$ $\mathbf{P}$-a.s. Then $\nu < \infty$ $\mathbf{P}^+$-a.s. and:*

(i) $(Q_1, Q_2, \ldots)$ *and* $(\widetilde{Q}_1, \widetilde{Q}_2, \ldots)$ *are identically distributed with respect to $\mathbf{P}^+$;*

(ii) $(\nu, Q_1, \ldots, Q_\nu)$ *and* $(\widetilde{Q}_1, \widetilde{Q}_2, \ldots)$ *are independent with respect to $\mathbf{P}^+$;*

(iii) $\mathbf{P}^+\{\nu = k, S_\nu \in dx\} = \mathbf{P}\{\iota = k, S_\iota \in dx\}$ *for all $k \geq 1$.*

PROOF. By monotonicity of $v$ and relation (2.16), we have

$$\mathbf{E}^+\left[\sum_{n=0}^\infty I_{\{S_n \leq x\}}\right] \leq v(x)\mathbf{E}\left[\sum_{n=0}^\infty I_{\{0 \leq S_n \leq x, L_n \geq 0\}}\right]$$

for every $x \geq 0$. The Markov property of the random walk $(S_n)_{n \geq 0}$ implies that the random sum on the right-hand side above has a geometrically decaying tail. Hence, its expectation is finite, which shows that under $\mathbf{P}^+$ the Markov chain $(S_n)_{n \geq 0}$ is transient. Consequently, $\mathbf{P}^+\{\nu < \infty\} = 1$.

To prove assertions (i) and (ii), we will first establish the corresponding statements for $S$. For $z \geq 0$, let

$$h_z(x) := \frac{v(x-z)}{v(x)}, \qquad x \geq 0.$$

Note that $0 \leq h_z \leq 1$, $h_z(x) = 0$ for $x < z$ and $h_z(x) \to 1$ as $x \to \infty$ [since the renewal function $v(x)$ satisfies $v(x) \to \infty$ and $v(x) - v(x-z) = O(1)$ as $x \to \infty$]. Moreover, from (1.5) we see that $h_z$ is harmonic with respect to the transition kernel $P^+$

$$\int P^+(x; dy) h_z(y) = h_z(x), \qquad x \geq z.$$

Since $\sigma_{z,k} := \min\{n \geq k : S_n < z\}, k \geq 0$ is a stopping time, the process $(h_z(S_{n \wedge \sigma_{z,k}}))_{n \geq k}$ is a martingale. Consequently, for $n \geq k$, we have

$$\mathbf{E}^+[h_z(S_{n \wedge \sigma_{z,k}})|S_k = x] = h_z(x).$$

Since $S_n \to \infty$ $\mathbf{P}^+$-a.s. as $n \to \infty$, the dominated convergence theorem entails

$$\mathbf{P}^+\{S_k, S_{k+1}, \ldots \geq z | S_k = x\} = \mathbf{E}^+\left[\lim_{n \to \infty} h_z(S_{n \wedge \sigma_{z,k}}) | S_k = x\right] = h_z(x).$$

It follows (with $x_0 = y_0 = 0$)

$$\mathbf{P}^+\{\nu = k, S_1 \in dx_1, \ldots, S_k \in dx_k, \widetilde{S}_1 \in dy_1, \ldots, \widetilde{S}_m \in dy_m\}$$

$$= \mathbf{1}_{\{x_1,\ldots,x_{k-1} > x_k\}} \mathbf{1}_{\{y_1,\ldots,y_m \geq 0\}} \prod_{i=1}^{k} P^+(x_{i-1}; dx_i)$$

$$\times \left(\prod_{j=1}^{m} P^+(y_{j-1} + x_k; dy_j + x_k)\right) h_{x_k}(y_m + x_k)$$

$$= \mathbf{1}_{\{x_1,\ldots,x_{k-1} > x_k\}} \left(\prod_{i=1}^{k} P^+(x_{i-1}; dx_i)\right) h_{x_k}(x_k) \prod_{j=1}^{m} P^+(y_{j-1}; dy_j)$$

$$= \mathbf{P}^+\{\nu = k, S_1 \in dx_1, \ldots, S_k \in dx_k\} \mathbf{P}^+\{S_1 \in dy_1, \ldots, S_m \in dy_m\}.$$

Thus, $(S_1, S_2, \ldots)$ and $(\widetilde{S}_1, \widetilde{S}_2, \ldots)$ are identical in distribution and $(\nu, S_1, \ldots, S_\nu)$ and $(\widetilde{S}_1, \widetilde{S}_2, \ldots)$ are independent (both with respect to $\mathbf{P}^+$).

Now we show that these properties carry over to the entire environment. By independence under the original measure, we have

(2.17) $\mathbf{P}\{Q_1 \in dq_1, \ldots, Q_k \in dq_k | S\} = k(X_1; dq_1) \cdots k(X_k; dq_k),$ $\mathbf{P}$-a.s.,

with $k(x; dq) := \mathbf{P}\{Q \in dq | X = x\}$. By definition of conditional expectation, using (2.16), we may conclude from (2.17) that

$$\mathbf{P}^+\{Q_1 \in dq_1, \ldots, Q_k \in dq_k | S\} = k(X_1; dq_1) \cdots k(X_k; dq_k), \qquad \mathbf{P}^+\text{-a.s.}$$

For Borel sets $B_i \subset \Delta$, the properties of $S$ established above imply

$$\mathbf{P}^+\{\nu = k, Q_1 \in B_1, \ldots, Q_k \in B_k, \widetilde{Q}_1 \in B_{k+1}, \ldots, \widetilde{Q}_m \in B_{k+m}\}$$

$$= \mathbf{E}^+\left[\prod_{i=1}^k k(X_i; B_i) \prod_{j=1}^m k(\widetilde{X}_j; B_{k+j}); \nu = k\right]$$

$$= \mathbf{E}^+\left[\prod_{i=1}^k k(X_i; B_i); \nu = k\right] \mathbf{E}^+\left[\prod_{j=1}^m k(X_j; B_{k+j})\right]$$

$$= \mathbf{P}^+\{\nu = k, Q_1 \in B_1, \ldots, Q_k \in B_k\} \mathbf{P}^+\{Q_1 \in B_{k+1}, \ldots, Q_m \in B_{k+m}\}.$$

Thus, we have proved (i) and (ii). As to (iii), using duality of random walks, we conclude

$$\begin{aligned}
\mathbf{P}^+&\{\nu = k, S_\nu \in dx\} \\
&= \mathbf{P}^+\{S_k \in dx, S_k - S_{k-1} < 0, \ldots, S_k - S_1 < 0\} h_x(x) \\
&= \mathbf{P}\{S_k \in dx, S_k - S_{k-1} < 0, \ldots, S_k - S_1 < 0\} \\
&= \mathbf{P}\{S_k \in dx, S_1 < 0, \ldots, S_{k-1} < 0\} \\
&= \mathbf{P}\{\iota = k, S_\iota \in dx\}.
\end{aligned}$$

This completes the proof of Lemma 2.6. □

2.3. *A convergent series theorem.* As an application of Tanaka's decomposition, we now prove a result, which previously had been obtained only under considerably stronger moment conditions [see Geiger and Kersting (2000), Kozlov (1976) and Vatutin and Dyakonova (2003)]. Let

$$\eta_k := \sum_{y=0}^\infty y(y-1) Q_k(\{y\}) \bigg/ \left(\sum_{y=0}^\infty y Q_k(\{y\})\right)^2, \qquad k \geq 1.$$

LEMMA 2.7. *Assume Assumptions* A1 *and* A2 *or* B1 *and* B2. *Then*

$$\sum_{k=0}^\infty \eta_{k+1} e^{-S_k} < \infty, \qquad \mathbf{P}^+\text{-}a.s.$$

PROOF. We will first estimate the $S_k$ from below. To this end, let $0 := \nu(0) < \nu(1) < \cdots$ be the times of prospective minima of $S$,

(2.18) $\quad \nu(j) := \min\{m > \nu(j-1) : S_{m+i} \geq S_m \text{ for all } i \geq 0\}, \qquad j \geq 1.$

Clearly,

(2.19) $\qquad\qquad\qquad S_k \geq S_{\nu(j)} \quad \text{if } k \geq \nu(j).$

By Lemma 2.6(i) and (ii), the random variable $S_{\nu(j)}$ is the sum of $j$ nonnegative i.i.d. random variables with positive mean. Thus, there exists some $c > 0$ such that

(2.20) $\qquad\qquad\qquad S_{\nu(j)} \geq cj \quad \text{eventually } \mathbf{P}^+\text{-a.s.}$

To get a lower bound on $\nu(j)$, observe that, by Lemma 2.6(i) and (ii), $\nu(j)$ is also the sum of $j$ nonnegative i.i.d. random variables, each with distribution $\nu = \nu(1)$. Lemma 2.6(iii) and (2.3) imply

$$\mathbf{P}^+\{\nu > k\} = \mathbf{P}\{\iota > k\} \leq \mathbf{P}\{\max(S_1, \ldots, S_k) \leq 0\} = o(k^{-\rho+\delta})$$

for every $\delta > 0$. Therefore, we have

$$\mathbf{E}^+ \nu^{\rho-\delta} < \infty \qquad \text{for all } \delta > 0.$$

Hence, an application of Theorem 13 in Chapter IX.3 in Petrov (1975) gives

$$(2.21) \qquad \nu(j) = O(j^{\rho^{-1}+\delta}), \qquad \mathbf{P}^+\text{-a.s.}$$

for every $\delta > 0$. Combining (2.19) and (2.20) with (2.21) gives

$$S_k \geq S_{\nu(\lfloor k^{\rho-\delta} \rfloor)} \geq c k^{\rho-\delta} \qquad \text{eventually } \mathbf{P}^+\text{-a.s.}$$

for every $\delta > 0$, which implies

$$(2.22) \qquad e^{-S_k} = O(e^{-k^{\rho-\delta}}), \qquad \mathbf{P}^+\text{-a.s.}$$

for all $\delta > 0$. To obtain this estimate, we have only used Assumption A1. Thus, it also holds under the stronger condition B1. However, under Assumption B1 it can be improved to

$$(2.23) \qquad e^{-S_k} = O(e^{-k^{\alpha^{-1}-\delta}}), \qquad \mathbf{P}^+\text{-a.s.}$$

for all $\delta > 0$. Recall from (1.6) that $\rho \leq \alpha^{-1}$. Hence, in view of (2.22), we may for the proof of (2.23) assume $0 < \alpha\rho < 1$. For this case Rogozin proved [cf. Chapter 8.9.2 in Bingham, Goldie and Teugels (1987), see also Doney (1995)] that the distribution of $S_\iota$ under $\mathbf{P}$ belongs to the domain of attraction of a stable law with index $\alpha\rho$. (We note that this holds for strict increasing ladder heights, as well as for weak increasing ladder heights, since the tails of both are identical up to a multiplicative constant.) Consequently, by Lemma 2.6(iii), for any $\delta > 0$, we have

$$\mathbf{P}^+\{S_\nu > x\} = \mathbf{P}\{S_\iota > x\} \geq x^{-\alpha\rho(1+\delta)},$$

if only $x$ is chosen large enough. Since $Y_i := S_{\nu(i)} - S_{\nu(i-1)}, i \geq 1$, are independent nonnegative random variables with the same distribution as $S_\nu$, we have

$$\mathbf{P}^+\{S_{\nu(j)} \leq j^{(1-\delta)/(\alpha\rho)}\}$$

$$\leq \mathbf{P}^+\left\{\max_{1 \leq i \leq j} Y_i \leq j^{(1-\delta)/(\alpha\rho)}\right\} = \mathbf{P}^+\{S_\nu \leq j^{(1-\delta)/(\alpha\rho)}\}^j$$

$$\leq \exp(-j\mathbf{P}^+\{S_\nu > j^{(1-\delta)/(\alpha\rho)}\}) \leq \exp(-j^{\delta^2}),$$

if only $j$ is large enough. The Borel–Cantelli lemma implies

$$S_{\nu(j)} \geq j^{(1-\delta)/(\alpha\rho)} \qquad \text{eventually } \mathbf{P}^+\text{-a.s.}$$

for all $\delta > 0$. Replacing (2.20) by this estimate assertion, (2.23) follows in much the same way as we derived (2.22).

The other part of the proof consists in estimating the $\eta_k$. First note that

$$\eta_k \leq \zeta_k(a) + \sum_{y=0}^{a-1} ayQ_k(\{y\}) \Big/ \left(\sum_{y=0}^{\infty} yQ_k(\{y\})\right)^2$$

$$\leq \zeta_k(a) + a\exp(-X_k)$$

for every $a \in \mathbb{N}_0$, where $\zeta_k(a)$ is the analogue of $\zeta(a)$ defined in terms of $Q_k$. Hence,

$$\sum_{k=0}^{\infty} \eta_{k+1}e^{-S_k} \leq \sum_{k=0}^{\infty} \zeta_{k+1}(a)e^{-S_k} + a\sum_{k=0}^{\infty} e^{-S_{k+1}}$$

(2.24)

$$\leq \sum_{k=0}^{\infty} (\zeta_{k+1}(a) + a)e^{-S_k},$$

and we are left with estimating the tail of $\zeta_k(a)$ under $\mathbf{P}^+$ for a suitable choice of $a$. Now note that the zero-delayed renewal function $v(x)$ satisfies the inequality $v(x+y) \leq v(x) + v(y)$. Therefore, by independence of the $Q_j$ under $\mathbf{P}$ and repeatedly using (2.16), we get

$$\mathbf{P}^+\{\zeta_k(a) > x\}$$

$$= \mathbf{E}[v(S_k); \zeta_k(a) > x, L_k \geq 0]$$

(2.25)
$$\leq \mathbf{E}[v(S_{k-1}) + v(X_k); \zeta_k(a) > x, L_{k-1} \geq 0]$$

$$= \mathbf{E}[v(S_{k-1}); L_{k-1} \geq 0]\mathbf{P}\{\zeta_k(a) > x\}$$

$$\quad + \mathbf{E}[v(X_k); \zeta_k(a) > x]\mathbf{P}\{L_{k-1} \geq 0\}$$

$$= \mathbf{P}\{\zeta(a) > x\} + \mathbf{E}[v(X); \zeta(a) > x]\mathbf{P}\{L_{k-1} \geq 0\}.$$

Now let $a \in \mathbb{N}_0$ and $\varepsilon > 0$ be such that Assumption A2 is satisfied. By means of (2.25) and the Markov inequality, it follows that, for every $x > 1$,

(2.26) $\quad \mathbf{P}^+\{\zeta_k(a) > x\} \leq \dfrac{c}{(\log x)^{(1/\rho)+\epsilon}} + \dfrac{c}{(\log x)^{1+\epsilon}}\mathbf{P}\{L_{k-1} \geq 0\}$

for some finite constant $c$. From the first part of Lemma 2.1 we see that

$$\mathbf{P}^+\{\zeta_k(a) > e^{k^{\rho-\delta'}}\}$$

$$= O(k^{-(\rho-\delta')(1/\rho+\epsilon)}) + O(k^{-(\rho-\delta')(1+\epsilon)}k^{-(1-\rho)+\delta'})$$

$$= O(k^{-1-\rho\epsilon/2}),$$

if only $\delta' > 0$ is chosen small enough. The Borel–Cantelli lemma implies

$$\zeta_k(a) = O(e^{k^{\rho-\delta'}}), \qquad \mathbf{P}^+\text{-a.s.}$$

for such $\delta'$. Combining this estimate with (2.22) and (2.24), the claim of the lemma follows under Assumptions A1 and A2.

Under Assumptions B1 and B2, the last estimate can again be further elaborated. Then, as $x \to \infty$,

$$v(x) = o(x^{\alpha(1-\rho)+\delta}) \tag{2.27}$$

for any $\delta > 0$. For the proof of (2.27) note that in analogy to (1.6) we have $\alpha(1-\rho) \leq 1$. Since in any case $v(x) = O(x)$, we may assume $0 < \alpha(1-\rho) < 1$. Then, by Rogozin's result, the distribution of $S_{\gamma_1}$ belongs to the domain of attraction of a stable law with index $\alpha(1-\rho)$. This implies that $v(x)$ is a regularly varying function with index $\alpha(1-\rho)$ [see Chapter 8.6.2 in Bingham, Goldie and Teugels (1987)], and (2.27) follows. Moreover, $\mathbf{E}|X|^{\alpha-\delta} < \infty$ for all $\delta > 0$, since by Assumption B1 the distribution of $X$ belongs to the domain of attraction of a stable law with index $\alpha$. Combining these two estimates gives

$$\mathbf{E}v(X)^{1/(1-\rho)-\delta} < \infty$$

for all $\delta > 0$. By means of Assumption B2 and Hölder's inequality, we obtain

$$\mathbf{E}[v(X)(\log^+ \zeta(a))^{\alpha\rho+\epsilon}] < \infty,$$

if only $\epsilon > 0$ is small enough and $a \in \mathbb{N}_0$ is sufficiently large. In view of Assumption B1, (2.25) and the Markov inequality, we get

$$\mathbf{P}^+\{\zeta_k(a) > x\} \leq \frac{c}{(\log x)^{\alpha+\epsilon}} + \frac{c}{(\log x)^{\alpha\rho+\epsilon}}\mathbf{P}\{L_{k-1} \geq 0\},$$

replacing the upper bound (2.26). Proceeding as above, we conclude

$$\begin{aligned}\mathbf{P}^+\{\zeta_k(a) > e^{k^{\alpha^{-1}-\delta'}}\} \\ = O(k^{-(1/\alpha-\delta')(\alpha+\epsilon)}) + O(k^{-(1/\alpha-\delta')(\alpha\rho+\epsilon)}k^{-(1-\rho)+\delta'}) \\ = O(k^{-1-\epsilon/(2\alpha)}),\end{aligned}$$

if $\delta' > 0$ is small enough. The Borel–Cantelli lemma implies

$$\zeta_k(a) = O(e^{k^{\alpha^{-1}-\delta'}}), \qquad \mathbf{P}^+\text{-a.s.}$$

for such $\delta'$. The claim of the lemma under Assumptions B1 and B2 follows from this estimate combined with (2.23) and (2.24). □

**3. Branching in conditioned environment.** Property (1.1) is unaffected under the change of measure, that is,

$$\mathbf{P}^+\{(Z_1,\ldots,Z_n) \in B|Z_0 = z_0, \Pi\} = \kappa_{n,z_0}(Q_1,\ldots,Q_n;B), \qquad \mathbf{P}^+\text{-a.s.}$$

This is an easy consequence of (1.2) and (2.16). Thus, $Z_0, Z_1, \ldots$ is still a branching process in a randomly fluctuating environment, however, the environment $Q_1, Q_2, \ldots$ is no longer built up from i.i.d. components. Let us call this a *branching process in conditioned environment*. Such processes exhibit a behavior, which is typical for supercritical branching processes. The following theorem states that, with respect to $\mathbf{P}^+$, the population has positive probability to survive forever. The statement holds for any initial distribution as long as $Z_0 \geq 1$ with positive probability.

PROPOSITION 3.1. *Assume Assumptions* A1 *and* A2 *or* B1 *and* B2. *Then*

$$\mathbf{P}^+\{Z_n > 0 \text{ for all } n | \Pi\} > 0, \qquad \mathbf{P}^+\text{-a.s.}$$

*In particular,*

$$\mathbf{P}^+\{Z_n > 0 \text{ for all } n\} > 0.$$

*Moreover, as* $n \to \infty$,

$$e^{-S_n} Z_n \to W^+, \qquad \mathbf{P}^+\text{-a.s.},$$

*where the random variable* $W^+$ *has the property*

$$\{W^+ > 0\} = \{Z_n > 0 \text{ for all } n\}, \qquad \mathbf{P}^+\text{-a.s.}$$

PROOF. In view of property (1.1), $Z_n$ is stochastically increasing with $Z_0$. Hence, for the proof of the first claim we may assume $Z_0 = 1$ $\mathbf{P}^+$-a.s. with no loss of generality. Consider the (random) generating functions

$$f_j(s) := \sum_{i=0}^{\infty} s^i Q_j(\{i\}), \qquad 0 \leq s \leq 1,$$

$j = 1, 2, \ldots$ and their compositions

(3.1) $$f_{k,n}(s) := f_{k+1}(f_{k+2}(\cdots f_n(s) \cdots)), \qquad 0 \leq k < n.$$

We note that the distributional identity (1.1) can be expressed as

(3.2) $$\mathbf{E}^+[s^{Z_n} | \Pi, Z_k] = f_{k,n}(s)^{Z_k}, \qquad \mathbf{P}^+\text{-a.s.}$$

We shall use an estimate on $f_{k,n}$ due to Agresti (1975) (see his Lemma 2), which was originally obtained through a comparison argument with linear fractional generating functions. A more direct proof may be given using the elementary identity

(3.3) $$\frac{1}{1 - f_{k,n}(s)} = \frac{e^{-(S_n - S_k)}}{1 - s} + \sum_{j=k}^{n-1} g_{j+1}(f_{j+1,n}(s)) e^{-(S_j - S_k)},$$

$$0 \leq s < 1,$$

with

$$g_j(s) := \frac{1}{1 - f_j(s)} - \frac{1}{f'_j(1)(1-s)}, \qquad 0 \leq s < 1,$$

and $g_j(1) := \lim_{s \to 1} g_j(s) = \eta_j/2$. Apparently, identity (3.3) has first been utilized by Jirina (1976). The coefficients possess the favorable property

$$0 \leq g_j(s) \leq \eta_j, \qquad 0 \leq s \leq 1,$$

which has been noticed by Geiger and Kersting (2000) (see their Lemma 2.6). Combining these formulas, we obtain Agresti's estimate

$$(3.4) \qquad f_{k,n}(s) \leq 1 - \left( \frac{e^{-(S_n - S_k)}}{1-s} + \sum_{j=k}^{n-1} \eta_{j+1} e^{-(S_j - S_k)} \right)^{-1}.$$

From (3.2) it follows that under the assumption $\mathbf{P}^+\{Z_0 = 1\} = 1$, we have

$$\mathbf{P}^+\{Z_n > 0 | \Pi\} = 1 - f_{0,n}(0), \qquad \mathbf{P}^+\text{-a.s.}$$

Recall that $S_n \to \infty$ $\mathbf{P}^+$-a.s. Hence, if we let $n \to \infty$, then (3.4) implies

$$\mathbf{P}^+\{Z_n > 0 \text{ for all } n | \Pi\} \geq \left( \sum_{j=0}^{\infty} \eta_{j+1} e^{-S_j} \right)^{-1}, \qquad \mathbf{P}^+\text{-a.s.}$$

Applying Lemma 2.7, we obtain

$$\mathbf{P}^+\{Z_n > 0 \text{ for all } n | \Pi\} > 0, \qquad \mathbf{P}^+\text{-a.s.}$$

This is the first claim of the proposition. The second claim is a well-known consequence of the martingale convergence theorem: Given the environment $\Pi$, $(Z_n/\mu_n)_{n \geq 0}$ is a martingale with respect to $\mathbf{P}^+$ and the filtration $\mathcal{F}$.

As to the proof of the third claim, note that $\mathbf{P}^+\{W^+ = 0\} \geq \mathbf{P}^+\{Z_n \to 0\}$, since $\{Z_n \to 0\} \subset \{W^+ = 0\}$. For the proof of the opposite inequality, we will use Tanaka's decomposition as an essential tool. To begin with, we show that

$$(3.5) \qquad \mathbf{P}^+\{Z_n \to 0 | \Pi\} + \mathbf{P}^+\{Z_n \to \infty | \Pi\} = 1, \qquad \mathbf{P}^+\text{-a.s.}$$

A sufficient condition for (3.5) is the following criterion [see Theorem 1 in Jagers (1974)]:

$$(3.6) \qquad \sum_{j=0}^{\infty} (1 - Q_j(\{1\})) = \infty, \qquad \mathbf{P}^+\text{-a.s.}$$

To verify (3.6), note that, by Lemma 2.6, the $Q_{\nu(k)+1}(\{1\})$, $k = 0, 1, \ldots$ are i.i.d. random variables. Also, since $\{Q_j(\{1\}) = 1\} \subset \{X_j = 0\}$ and the case of a degenerate random walk $S$ is excluded by Spitzer's condition A1, we have

$$\mathbf{P}^+\{Q_{\nu(0)+1}(\{1\}) = 1\} \leq \mathbf{P}^+\{X_1 = 0\} < 1.$$

Hence,
$$\sum_{j=1}^{\infty}(1-Q_j(\{1\})) \geq \sum_{k=0}^{\infty}(1-Q_{\nu(k)+1}(\{1\})) = \infty, \qquad \mathbf{P}^+\text{-a.s.}$$

Clearly, (3.5) implies

(3.7) $$\mathbf{P}^+\{Z_n \to 0\} + \mathbf{P}^+\{Z_n \to \infty\} = 1.$$

Now observe that from (3.2) and (3.4) we get

(3.8)
$$\mathbf{E}^+[\exp(-\lambda e^{-S_n}Z_n)|Z_k=1,\Pi]$$
$$= f_{k,n}(\exp(-\lambda e^{-S_n}))$$
$$\leq 1 - \left(\frac{e^{-(S_n-S_k)}}{1-\exp(-\lambda e^{-S_n})} + \sum_{j=k}^{n-1}\eta_{j+1}e^{-(S_j-S_k)}\right)^{-1}, \qquad \mathbf{P}^+\text{-a.s.}$$

for every $\lambda \geq 0$ and $k < n$. Recall that $S_n \to \infty$ and $e^{-S_n}Z_n \to W^+$ $\mathbf{P}^+$-a.s. Hence, letting first $n \to \infty$ and then $\lambda \to \infty$ gives

$$\mathbf{P}^+\{W^+=0|Z_k=1,\Pi\} \leq 1 - \left(\sum_{j=k}^{\infty}\eta_{j+1}e^{-(S_j-S_k)}\right)^{-1}, \qquad \mathbf{P}^+\text{-a.s.}$$

Since the times $\nu(k)$ of prospective minima are determined by the environment only, we may replace $k$ by $\nu(k)$ in the last estimate. Moreover, identity (1.1) implies

$$\mathbf{P}^+\{W^+=0|Z_{\nu(k)}=j,\Pi\} = \mathbf{P}^+\{W^+=0|Z_{\nu(k)}=1,\Pi\}^j.$$

Combining these observations gives

$$\mathbf{P}^+\{W^+=0|\Pi\}$$
$$= \mathbf{E}^+[\mathbf{P}^+\{W^+=0|Z_{\nu(k)},\Pi\}|\Pi]$$
$$\leq \mathbf{E}^+\left[\left(1-\frac{1}{\sum_{j=\nu(k)}^{\infty}\eta_{j+1}e^{-(S_j-S_{\nu(k)})}}\right)^{Z_{\nu(k)}}\bigg|\Pi\right]$$
$$\leq \mathbf{P}^+\{Z_{\nu(k)} \leq z|\Pi\} + \left(1-\frac{1}{\sum_{j=\nu(k)}^{\infty}\eta_{j+1}e^{-(S_j-S_{\nu(k)})}}\right)^z, \qquad \mathbf{P}^+\text{-a.s.}$$

for every $z \geq 0$. By Lemma 2.6, the law of the second term on the right-hand side above does not depend on $k$. Hence, taking first expectations and then letting $k \to \infty$, we see from (3.7) that

$$\mathbf{P}^+\{W^+=0\} \leq \mathbf{P}^+\{Z_n \to 0\} + \mathbf{E}^+\left(1-\frac{1}{\sum_{j=0}^{\infty}\eta_{j+1}e^{-S_j}}\right)^z.$$

Finally, letting $z \to \infty$, an application of Lemma 2.7 yields

$$\mathbf{P}^+\{W^+ = 0\} - \mathbf{P}^+\{Z_n \to 0\} \leq \mathbf{P}^+\left(\sum_{j=0}^{\infty} \eta_{j+1} e^{-S_j} = \infty\right) = 0,$$

which completes the proof of Proposition 3.1.  □

**4. Proofs of Theorems 1.1 and 1.3–1.5.** The general approach of our proofs is to replace the conditioning event $\{Z_n > 0\}$ by other events, which are easier to handle. This strategy has been used before: Kozlov (1976) considered the event that only a few descending ladder epochs of the random walk $S$ occur before time $n$ and Geiger and Kersting (2000) conditioned on the event that the random walk has a high minimum $L_n$. We follow the approach of Dyakonova, Geiger and Vatutin (2004) and condition on the event that $S$ attains its minimal value extraordinarily early, which is conceptually more appealing and also allows some simplifications in the proofs. The next lemma presents a main argument, which will be used throughout the proofs of our theorems. Recall the definitions of $\tau_n$ and $L_{k,n}$ from (1.9) and (2.8).

LEMMA 4.1.  *Assume Assumption* A1 *and let* $m \in \mathbb{N}_0$. *Suppose* $V_1, V_2, \ldots$ *is a uniformly bounded sequence of real-valued random variables, which, for every* $k \geq 0$, *satisfy*

$$(4.1)\ \mathbf{E}[V_n; Z_{k+m} > 0, L_{k,n} \geq 0 | \mathcal{F}_k] = \mathbf{P}\{L_n \geq 0\}(V_{k,\infty} + o(1)), \qquad \mathbf{P}\text{-}a.s.$$

*with random variables* $V_{1,\infty} = V_{1,\infty}(m)$, $V_{2,\infty} = V_{2,\infty}(m), \ldots$. *Then*

$$(4.2) \quad \mathbf{E}[V_n; Z_{\tau_n+m} > 0] = \mathbf{P}\{L_n \geq 0\}\left(\sum_{k=0}^{\infty} \mathbf{E}[V_{k,\infty}; \tau_k = k] + o(1)\right),$$

*where the right-hand side series is absolutely convergent.*

In our applications of Lemma 4.1, the $V_n$ will be typically of the form $V_n = U_n I_{\{Z_n > 0\}}$ with random $U_n$. Relation (4.2) then reflects the fact that, given survival at generation $n$, the history of the branching process splits into two independent pieces. The summands display the evolution of the branching process up to time $k$ when $\tau_n = k$, whereas the common factor $\mathbf{P}\{L_n \geq 0\}$ arises from the evolution after time $\tau_n$.

PROOF OF LEMMA 4.1.  Fix $m \in \mathbb{N}_0$. We may assume $0 \leq V_n \leq 1$ since assumption (4.1) implies the corresponding statements for the positive and the negative part of the $V_n$. Using first (2.9) and then the independence of the $X_j$ and the estimate (1.3), we obtain

$$\mathbf{E}[V_n; Z_{\tau_n+m} > 0, \tau_n > l] \leq \mathbf{P}\{Z_{\tau_n} > 0, \tau_n > l\}$$
$$= \sum_{k=l+1}^{n} \mathbf{P}\{Z_k > 0, \tau_k = k, L_{k,n} \geq 0\}$$
$$\leq \sum_{k=l+1}^{n} \mathbf{E}[e^{S_k}; \tau_k = k]\mathbf{P}\{L_{n-k} \geq 0\}$$

for every $l \in \mathbb{N}_0$. Applying Lemma 2.2 with $u(x) := e^{-x}$ gives

(4.3)    $\lim_{l \to \infty} \limsup_{n \to \infty} (\mathbf{P}\{L_n \geq 0\})^{-1} \mathbf{E}[V_n; Z_{\tau_n+m} > 0, \tau_n > l] = 0.$

On the other hand, using (2.9) again, we have

(4.4)
$$\mathbf{E}[V_n; Z_{\tau_n+m} > 0, \tau_n = k]$$
$$= \mathbf{E}[V_n; Z_{k+m} > 0, \tau_k = k, L_{k,n} \geq 0]$$
$$= \mathbf{E}[\mathbf{E}[V_n; Z_{k+m} > 0, L_{k,n} \geq 0|\mathcal{F}_k]; \tau_k = k]$$

for every $k \leq n$. Now observe that, by independence of the $X_j$, we get

$$\mathbf{E}[V_n; Z_{k+m} > 0, L_{k,n} \geq 0|\mathcal{F}_k]$$
$$\leq \mathbf{P}\{L_{k,n} \geq 0|\mathcal{F}_k\} = \mathbf{P}\{L_{n-k} \geq 0\}, \qquad \mathbf{P}\text{-a.s.}$$

Since $\mathbf{P}\{L_{n-k} \geq 0\} \sim \mathbf{P}\{L_n \geq 0\}$ for fixed $k$, relation (4.4) and the dominated convergence theorem, combined with the assumption of the lemma, imply

(4.5)    $\lim_{n \to \infty} (\mathbf{P}\{L_n \geq 0\})^{-1} \mathbf{E}[V_n; Z_{\tau_n+m} > 0, \tau_n = k] = \mathbf{E}[V_{k,\infty}; \tau_k = k]$

for every $k \in \mathbb{N}_0$. Consequently,

(4.6)    $\sum_{k=l+1}^{\infty} \mathbf{E}[V_{k,\infty}; \tau_k = k] \leq \limsup_{n \to \infty} (\mathbf{P}\{L_n \geq 0\})^{-1} \mathbf{E}[V_n; Z_{\tau_n+m} > 0, \tau_n > l]$

for every $l \in \mathbb{N}_0$. By means of the triangle inequality, we obtain from (4.5) and (4.6)

(4.7)
$$\limsup_{n \to \infty} \left| (\mathbf{P}\{L_n \geq 0\})^{-1} \mathbf{E}[V_n; Z_{\tau_n+m} > 0] - \sum_{k=0}^{\infty} \mathbf{E}[V_{k,\infty}; \tau_k = k] \right|$$
$$\leq 2 \limsup_{n \to \infty} (\mathbf{P}\{L_n \geq 0\})^{-1} \mathbf{E}[V_n; Z_{\tau_n+m} > 0, \tau_n > l]$$

for every $l \in \mathbb{N}_0$. Since the left-hand side of (4.7) does not depend on $l$, the claim of the lemma follows from (4.3) by letting $l \to \infty$ in (4.7). □

For convenience we introduce the notation

$$A_{\text{u.s.}} := \{Z_n > 0 \text{ for all } n \geq 0\}$$

for the event of ultimate survival.

PROOF OF THEOREM 1.1. For $z, n \in \mathbb{N}_0$, we write

$$\psi(z, n) := \mathbf{P}_z\{Z_n > 0, L_n \geq 0\}.$$

Note that $\psi(0, n) = 0$. Choosing $Y_n = I_{\{Z_n > 0\}}$ and $Y_\infty = I_{A_{\text{u.s.}}}$ in Lemma 2.5, we get, for $z \geq 1$,

$$(4.8) \qquad \psi(z, n) \sim \mathbf{P}\{L_n \geq 0\} \mathbf{P}_z^+\{A_{\text{u.s.}}\} \qquad \text{as } n \to \infty.$$

Furthermore, for $k \leq n$, we have

$$(4.9) \qquad \mathbf{P}\{Z_n > 0, L_{k,n} \geq 0 | \mathcal{F}_k\} = \psi(Z_k, n - k), \qquad \mathbf{P}\text{-a.s.}$$

Relations (4.8) and (4.9) show that we may apply Lemma 4.1 to $V_n = I_{\{Z_n > 0\}}$, $V_{k,\infty} = \mathbf{P}_{Z_k}^+\{A_{\text{u.s.}}\}$ and $m = 0$ to obtain

$$\mathbf{P}\{Z_n > 0\} \sim \theta \mathbf{P}\{L_n \geq 0\} \qquad \text{as } n \to \infty,$$

where

$$(4.10) \qquad \theta := \sum_{k=0}^{\infty} \mathbf{E}[\mathbf{P}_{Z_k}^+\{A_{\text{u.s.}}\}; \tau_k = k] < \infty.$$

For $\theta$ being strictly positive, note that Proposition 3.1 implies $\mathbf{P}_z^+\{A_{\text{u.s.}}\} > 0$ for all $z \geq 1$. □

REMARK. It is interesting to note that the sum in representation (4.10) of $\theta$ can be interpreted as follows: Call a strict descending ladder epoch of the associated random walk an unfavorable generation (at such epochs the probability of survival is particularly low). If the members of each unfavorable generation are transfered into a conditioned random environment and branch according to this new environment, then $\theta$ is the expected number of such clans, which survive forever.

PROOF OF THEOREM 1.3. Let $\phi$ be a bounded continuous function on the space $D[0, 1]$ of càdlàg functions on the unit interval. For $s \in \mathbb{R}$, let $W^s$ denote the process with constant paths

$$W_t^s := e^{-s} W^+, \qquad 0 \leq t \leq 1,$$

where $W^+$ is specified in Proposition 3.1. For fixed $s \in \mathbb{R}$, Proposition 3.1 shows that, as $n, r_n \to \infty$ with $r_n \leq n$, the process $e^{-s} X^{r_n, n}$ converges to $W^s$ in the metric of uniform convergence and, consequently, in the Skorokhod-metric on the space $D[0,1]$ $\mathbf{P}^+$-a.s.,

$$Y_n := \phi(e^{-s}X^{r_n,n})I_{\{Z_n>0\}} \to Y_\infty := \phi(W^s)I_{\{W^+>0\}}, \qquad \mathbf{P}^+\text{-a.s.}$$

(In fact, since the limiting process $W^s$ has continuous paths, convergence in the two metrics is equivalent.) For $r \leq n$ and $z \in \mathbb{N}_0$, define

$$\psi(z,s,r,n) := \mathbf{E}_z[\phi(e^{-s}X^{r,n}); Z_n > 0, L_n \geq 0].$$

Lemma 2.5 entails

$$\psi(z,s,r_n,n) = \mathbf{P}\{L_n \geq 0\}(\mathbf{E}_z^+[\phi(W^s); W^+ > 0] + o(1)).$$

Now observe that, for $k \leq r \leq n$,

$$\mathbf{E}[\phi(X^{r,n}); Z_n > 0, L_{k,n} \geq 0 | \mathcal{F}_k] = \psi(Z_k, S_k, r-k, n-k), \qquad \mathbf{P}\text{-a.s.}$$

Thus, we may apply Lemma 4.1 to the random variables $V_n = \phi(X^{r_n,n})I_{\{Z_n>0\}}$ and $V_{k,\infty} = \mathbf{E}_{Z_k}^+[\phi(W^{S_k}); W^+ > 0]$ with $m = 0$. Also using Theorem 1.1, we obtain

$$\mathbf{E}[\phi(X^{r_n,n})|Z_n > 0] \to \int \phi(w)\lambda(dw) \qquad \text{as } n \to \infty,$$

where $\lambda$ is the measure on the space of càdlàg functions on $[0,1]$ given by

$$\lambda(dw) := \frac{1}{\theta}\sum_{k=0}^\infty \mathbf{E}[\lambda_{Z_k,S_k}(dw); Z_k > 0, \tau_k = k]$$

with

$$\lambda_{z,s}(dw) := \mathbf{P}_z^+[W^s \in dw, W^+ > 0].$$

By Proposition 3.1, the total mass of $\lambda_{z,s}$ is $\mathbf{P}^+\{A_{\text{u.s.}}\}$. Hence, the representation of $\theta$ in (4.10) shows that $\lambda$ is a probability measure. Again using Proposition 3.1, we see that $\lambda_{z,s}$ puts its entire mass on strictly positive constant functions and, hence, so does $\lambda$. This completes the proof of the theorem. □

PROOF OF THEOREM 1.4. Note that $\min(S_0, \ldots, S_n) = L_n \wedge 0$. We consider $V_n = \phi(\tau_n, L_n \wedge 0)I_{\{Z_n>0\}}$ for some bounded measurable function $\phi$ on $\mathbb{N}_0 \times \mathbb{R}_0^-$. Since $L_{k,n} \geq 0$ implies $\tau_k = \tau_n$ [cf. (2.9)], we have

$$\mathbf{E}[V_n; Z_k > 0, L_{k,n} \geq 0 | \mathcal{F}_k] = \mathbf{E}[\phi(\tau_k, S_{\tau_k}); Z_n > 0, L_{k,n} \geq 0 | \mathcal{F}_k]$$
$$= \phi(\tau_k, S_{\tau_k})\mathbf{P}\{Z_n > 0, L_{k,n} \geq 0 | \mathcal{F}_k\}, \qquad \mathbf{P}\text{-a.s.}$$

Thus, we may apply Lemma 4.1 in just the same manner as in the proof of Theorem 1.1 with $V_{k,\infty} = \phi(\tau_k, S_{\tau_k})\mathbf{P}^+_{Z_k}\{A_{\text{u.s.}}\}$ and obtain

$$\mathbf{E}[\phi(\tau_n, L_n \wedge 0)|Z_n > 0]$$
$$\to \frac{1}{\theta}\sum_{k=0}^{\infty}\mathbf{E}[\phi(k, S_k)\mathbf{P}^+_{Z_k}\{A_{\text{u.s.}}\}; \tau_k = k] \qquad \text{as } n \to \infty.$$

This entails the desired result.  □

PROOF OF THEOREM 1.5.  Let $k, m \geq 0$ and $k + m \leq n$. Similar to the proof of Lemma 2.3, we may, in view of (2.7), decompose the stochastic process under consideration as $S^n = S^{k+m,n} + \widetilde{S}^{k+m,n}$. Let $\phi$ be a bounded continuous function on $D[0,1]$ and define

$$\psi(w, x) := \mathbf{E}[\phi(w + \widetilde{S}^{k+m,n}); L_{k+m,n} \geq -x]$$

for $w \in D[0,1]$ and $x \geq 0$. By Lemma 2.3, given $L_{k+m,n} \geq -x$, the process $\widetilde{S}^{k+m,n}$ converges in distribution to $L^+$ as $n \to \infty$ for each $k$ and $m$. Hence, if the càdlàg functions $w^n$ converge uniformly to the zero function, then

$$\psi(w^n, x) = \mathbf{P}\{L_{n-(k+m)} \geq -x\}(\mathbb{E}\phi(L^+) + o(1))$$
$$= v(x)\mathbf{P}\{L_n \geq 0\}(\mathbb{E}\phi(L^+) + o(1)),$$

where for the second equality we have used (2.1). Since

$$(4.11) \qquad \{L_{k,n} \geq 0\} = \{L_{k,k+m} \geq 0\} \cap \{L_{k+m,n} \geq -(S_{k+m} - S_k)\}$$

and since $S^{k+m,n}$ converges uniformly to zero as $n \to \infty$ **P**-a.s., we obtain

$$\mathbf{E}[\phi(S^n); Z_{k+m} > 0, L_{k,n} \geq 0|\mathcal{F}_{k+m}]$$
$$(4.12) \qquad = \psi(S^{k+m,n}, S_{k+m} - S_k)I_{\{Z_{k+m}>0, L_{k,k+m} \geq 0\}}$$
$$= v(S_{k+m} - S_k)\mathbf{P}\{L_n \geq 0\}$$
$$\times (\mathbb{E}\phi(L^+) + o(1))I_{\{Z_{k+m}>0, L_{k,k+m}\geq 0\}}, \qquad \mathbf{P}\text{-a.s.}$$

From (2.2) and (4.11) we deduce

$$|\mathbf{E}[\phi(S^n); Z_{k+m} > 0, L_{k,n} \geq 0|\mathcal{F}_{k+m}]|$$
$$\leq \sup |\phi|\mathbf{P}\{L_{k,n} \geq 0|\mathcal{F}_{k+m}\}$$
$$= \sup |\phi|\mathbf{P}\{L_{k+m,n} \geq -(S_{k+m} - S_k)|\mathcal{F}_{k+m}\}I_{\{L_{k,k+m}\geq 0\}}$$
$$\leq cv(S_{k+m} - S_k)\mathbf{P}\{L_{n-(k+m)} \geq 0\}I_{\{L_{k,k+m}\geq 0\}}, \qquad \mathbf{P}\text{-a.s.}$$

for some $c > 0$. Also, $\mathbf{E}[v(S_{k+m} - S_k); L_{k,k+m} \geq 0|\mathcal{F}_k] = v(0) < \infty$ **P**-a.s., by (1.5). Hence, by means of the dominated convergence theorem and (2.16),

we conclude from (4.12)

$$\mathbf{E}[\phi(S^n); Z_{k+m} > 0, L_{k,n} \geq 0 | \mathcal{F}_k]$$
$$= (\mathbb{E}\phi(L^+) + o(1))\mathbf{P}\{L_n \geq 0\}$$
$$\times \mathbf{E}[v(S_{k+m} - S_k); Z_{k+m} > 0, L_{k,k+m} \geq 0 | \mathcal{F}_k]$$
$$= (\mathbb{E}\phi(L^+) + o(1))\mathbf{P}\{L_n \geq 0\}\mathbf{P}^+_{Z_k}\{Z_m > 0\}, \qquad \mathbf{P}\text{-a.s.}$$

Applying Lemma 4.1 to $V_n = \phi(S^n)$ gives

$$\mathbf{E}[\phi(S^n); Z_{\tau_n+m} > 0] = (\mathbb{E}\phi(L^+) + o(1))\mathbf{P}\{L_n \geq 0\}\sum_{k=0}^{\infty}\mathbf{E}[\mathbf{P}^+_{Z_k}\{Z_m > 0\}; \tau_k = k].$$

In particular, we have

$$(4.13) \quad \mathbf{P}\{Z_{\tau_n+m} > 0\} \sim \mathbf{P}\{L_n \geq 0\}\sum_{k=0}^{\infty}\mathbf{E}[\mathbf{P}^+_{Z_k}\{Z_m > 0\}; \tau_k = k],$$

where the right-hand side series is convergent. Now observe that

$$|\mathbb{E}\phi(L^+)\mathbf{P}\{Z_n > 0\} - \mathbf{E}[\phi(S^n); Z_n > 0]|$$
$$\leq |\mathbb{E}\phi(L^+)\mathbf{P}\{Z_n > 0\} - \mathbf{E}[\phi(S^n); Z_{\tau_n+m} > 0]|$$
$$+ \sup|\phi|\mathbb{E}|I_{\{Z_n>0\}} - I_{\{Z_{\tau_n+m}>0\}}|$$

and

$$\mathbb{E}|I_{\{Z_n>0\}} - I_{\{Z_{\tau_n+m}>0\}}|$$
$$\leq (\mathbf{P}\{Z_n > 0\} - \mathbf{P}\{Z_{n+m} > 0\}) + (\mathbf{P}\{Z_{\tau_n+m} > 0\} - \mathbf{P}\{Z_{n+m} > 0\}).$$

Combining these estimates with Theorem 1.1 gives

$$(4.14) \quad \begin{aligned}&|\mathbb{E}\phi(L^+) - \mathbf{E}[\phi(S^n)|Z_n > 0]|\\&\leq 2\sup|\phi|\left(\frac{1}{\theta}\sum_{k=0}^{\infty}\mathbf{E}[\mathbf{P}^+_{Z_k}\{Z_m > 0\}; \tau_k = k] - 1\right) + o(1).\end{aligned}$$

By the dominated convergence theorem and (4.10),

$$\sum_{k=0}^{\infty}\mathbf{E}[\mathbf{P}^+_{Z_k}\{Z_m > 0\}; \tau_k = k] \downarrow \theta \qquad \text{as } m \to \infty.$$

Since the left-hand side of (4.14) does not depend on $m$, the assertion of Theorem 1.5 follows. $\square$

V. I. Afanasyev
V. A. Vatutin
Department of Discrete Mathematics
Steklov Institute
8 Gubkin Street
117 966 Moscow
GSP-1
Russia
e-mail: iafan@mail.ru
e-mail: vatutin@mi.ras.ru

J. Geiger
Fachbereich Mathematik
Technische Universität Kaiserslautern
Postfach 3049
D-67653 Kaiserslautern
Germany
e-mail: jgeiger@mathematik.uni-kl.de

G. Kersting
Fachbereich Mathematik
Universität Frankfurt
Fach 187
D-60054 Frankfurt am Main
Germany
e-mail: kersting@math.uni-frankfurt.de